\documentstyle[epsf]{article}

\newcommand{\bigzerou}{%
\smash{\lower1.7ex\hbox{\bg 0}}}
\newtheorem{theorem}{Theorem}
\newtheorem{prop}{Proposition}

\newcommand{\ba}{\begin{eqnarray}}
\newcommand{\ea}{\end{eqnarray}}
\newcommand{\no}{\nonumber}

\newcommand{\mapright}[1]{%
\smash{\mathop{%
\hbox to 1.0cm{\rightarrowfill}}\limits^{#1}}}
\newcommand{\mapleft}[1]{%
\smash{\mathop{%
\hbox to 1.3cm{\leftarrowfill}}\limits^{#1}}}

\begin{document}
\title{
\begin{flushright}
  \begin{minipage}[b]{5em}
    \normalsize
    ${}$      \\
  \end{minipage}
\end{flushright}
{\bf Classical Computation of Number of Lines in Projective Hypersurfaces: Origin of Mirror Transformation}}
\author{Masao Jinzenji\\
\\
\it Division of Mathematics, Graduate School of Science \\
\it Hokkaido University \\
\it  Kita-ku, Sapporo, 060-0810, Japan\\
{\it e-mail address: jin@math.sci.hokudai.ac.jp}}
\maketitle
\begin{abstract}
In this paper, we discuss classical derivation of the residue integral representation of the $d=1$ rational Gromov-Witten 
invariants of projective hypersurfaces that followed from localization technique.
\end{abstract}
\section{Introduction}
Geometrical origin of the mirror transformation in the mirror computation of the Gromov-Witten invariants is considered as
an effect of the difference between the moduli space of Gauged linear Sigma Model and the one of Non-Linear Sigma Model.
In this paper, we focus on this problem in the case of the $d=1$ rational Gromov-Witten invariants of projective hypersurface,
which is the simplest example of the effect of the mirror transformation.

In \cite{mmg}, we introduced the compactified moduli space of polynomial maps with two marled points $\widetilde{Mp}_{0,2}(N,d)$, 
which is expected to play the role of the moduli space of Gauged Linear Sigma Model whose target space is $CP^{N-1}$. In the $d=1$ 
case, $\widetilde{Mp}_{0,2}(N,1)$ turns out to be $CP^{N-1}\times CP^{N-1}$. We then introduced the two point function $w({\cal O}_{h^a}
{\cal O}_{h_b})_{0,1}$ which is considered as the correlation function of the Gauged Linear Sigma Model whose target space is the 
degree $k$ hypersurface in $CP^{N-1}$ (we denote this hypersurface by $M_{N}^{k}$). 
Here $h$ is the hyperplane class of $CP^{N-1}$. Explicitly, it is given by the following 
formula:
\begin{equation}
w({\cal O}_{h^a}{\cal O}_{h^b})_{0,1}=\int_{CP^{N-1}\times CP^{N-1}}(h_1)^{a}\cdot\bigl(\prod_{j=0}^{k}(jh_1+(k-j)h_2)\bigr)\cdot(h_{2})^b,
\label{w}
\end{equation}
where $h_1$ (resp. $h_2$) is the hyperplane class of the first (resp. second) $CP^{N-1}$. The r.h.s. of (\ref{w}) can be written as the 
residue integral:
\begin{equation}
w({\cal O}_{h^a}{\cal O}_{h^b})_{0,1}=\frac{1}{(2\pi\sqrt{-1})^2}\oint_{C_0}\frac{dz_1}{(z_1)^N}
\oint_{C_0}\frac{dz_2}{(z_2)^N}(z_1)^{a}\cdot\bigl(\prod_{j=0}^{k}(jz_1+(k-j)z_2)\bigr)\cdot(z_{2})^b,
\label{wr}
\end{equation}    
where $\displaystyle{\frac{1}{2\pi\sqrt{-1}}\oint_{C_0}dz}$ represents the operation of taking a residue at $z=0$.
On the other hand, we derived in \cite{dpr} the residue integral representation of the two point rational Gromof-Witten invariant 
$\langle{\cal O}_{h^a}{\cal O}_{h^b}\rangle_{0,1}$
of $M_{N}^{k}$ by taking the non-equivariant limit of the result of localization technique by Kontsevich \cite{kont}:
\begin{equation}
\langle{\cal O}_{h^a}{\cal O}_{h^b}\rangle_{0,1}=-\frac{1}{2}\cdot\frac{1}{(2\pi\sqrt{-1})^2}\oint_{C_0}\frac{dz_1}{(z_1)^N}
\oint_{C_0}\frac{dz_2}{(z_2)^N}\bigl(\prod_{j=0}^{k}(jz_1+(k-j)z_2)\bigr)\cdot(z_1-z_2)^2\cdot\frac{(z_1)^a-(z_2)^a}{z_1-z_2}\cdot
\frac{(z_1)^b-(z_2)^b}{z_1-z_2}.
\label{gr}
\end{equation} 
From (\ref{wr}) and (\ref{gr}), we can derive the equality:
\begin{equation}
\langle{\cal O}_{h^a}{\cal O}_{h^b}\rangle_{0,1}=w({\cal O}_{h^a}{\cal O}_{h^b})_{0,1}-w({\cal O}_{h^{a+b}}{\cal O}_{1})_{0,1},
\end{equation}
which is the simplest example of the mirror transformation.

In this paper, we discuss another proof of (\ref{gr}) by using the classical cohomology theory of the Grassmann variety $G(2,N)$,
which parametrizes $2$-dimensional linear subspaces in ${\bf C}^{N}$.            
Let us introduce the following rational functions in $z$ and $w$ to state our result in this paper: 
\begin{eqnarray}
e^{k}(z,w)&:=&\prod_{j=0}^{k}(jz+(k-j)w),\no\\
w_a(z,w)&:=&\frac{z^a-w^a}{z-w}.
\end{eqnarray}
With this set-up, we state the following theorem on the $d=1$ rational Gromov-Witten invariants of $M_{N}^{k}$. 
\begin{theorem}
\begin{equation}
\langle\prod_{j=1}^{n}{\cal O}_{h^{a_j}}\rangle_{0,1}=\frac{1}{(2\pi\sqrt{-1})^2}\oint_{C_0}\frac{dz_1}{(z_1)^N}
\oint_{C_0}\frac{dz_2}{(z_2)^N}(z_2-z_1)\cdot
e^{k}(z_1,z_2)
\cdot(z_1)^{a_1}\prod_{j=2}^{n}w_{a_j}(z_1,z_2),
\label{thm}
\end{equation}
where $\displaystyle{\frac{1}{2\pi\sqrt{-1}}\oint_{C_0}dz}$ represents the operation of taking a residue at $z=0$.
\end{theorem}
If $n=2$, (\ref{thm}) reduces to (\ref{gr}). To prove Theorem 1, we use the formula that represents $\langle\prod_{j=1}^{n}{\cal O}_{h^{a_j}}\rangle_{0,1}$ in terms of cohomology classes of $G(2,N)$:
\begin{equation}
\langle\prod_{j=1}^{n}{\cal O}_{h^{a_j}}\rangle_{0,1}=
\int_{G(2,N)}c_{top}(S^k(S_G^{*}))\cdot(\prod_{j=1}^{n}\sigma^{*}_{a_j-1}\bigr),
\label{ingw}
\end{equation}
where $S_G$ is the tautological vector bundle of $G(2,N)$ and $\sigma_{a}^{*}$ is the Schubert class that will be introduced in 
Section 3. Then we consider the following sequence of maps:
\begin{equation}
(CP^{N-1}\times CP^{N-1}-\Delta)\stackrel{f}{\longrightarrow}{\bf P}({\bf C}^{N}/S_1)\stackrel{\eta}{\longrightarrow} G(2,N)
\label{map}
\end{equation}
where $\Delta$ is the diagonal set of $CP^{N-1}\times CP^{N-1}$, $S_1$ is the tautological line bundle of $CP^{N-1}$ and 
${\bf P}({\bf C}^{N}/S_1)$ is the projective bundle over $CP^{N-1}$ that will be introduced in Section 2. Let 
$U$ be a set of pairs of two linearly independent vectors in ${\bf C}^{N}$:
\begin{equation}
U=\{({\bf a}_1,{\bf a}_2)\;|\; {\bf a}_1,{\bf a}_2\in {\bf C}^{N},\;\mbox{${\bf a}_1$ and ${\bf a}_2$ are linearly independent}\}.
\end{equation}  
Since $\left[\begin{array}{cc}\alpha&\beta\\ \gamma&\delta\end{array}\right]\in GL(2,C)$ acts on $U$ by 
$({\bf a}_1,{\bf a}_2)\left[\begin{array}{cc}\alpha&\beta\\ \gamma&\delta\end{array}\right]=(\alpha{\bf a}_1+\gamma{\bf a}_2,\beta{\bf a}_1+\delta{\bf a}_2)$, we have,
\begin{equation}
(CP^{N-1}\times CP^{N-1}-\Delta)=U/B_1,\;{\bf P}({\bf C}^{N}/S_1)=U/B_2,\;G(2,N)=U/GL(2,C),
\end{equation}
where $B_1$ and $B_2$ are the subgroups of $GL(2,C)$ given by,
\begin{eqnarray}
&& B_1=\{\left[\begin{array}{cc}\alpha&0\\ 0&\delta\end{array}\right]\;\;|\;\;\alpha\delta\neq 0\;\;\},\no\\
&& B_2=\{\left[\begin{array}{cc}\alpha&\beta\\ 0&\delta\end{array}\right]\;\;|\;\;\alpha\delta\neq 0\;\;\}.
\end{eqnarray} 
Therefore, $f$ and $\eta$  in (\ref{map}) are natural projections. The main idea to prove Theorem 1 is to rewrite 
the r.h.s. of (\ref{ingw}) as the intersection number of ${\bf P}({\bf C}^{N}/S_1)$ via the projection formula with respect 
to $\eta$. But the final result can be interpreted as the intersection number of $CP^{N-1}\times CP^{N-1}$ because 
the residue integral representation indicates that the r.h.s. of (\ref{thm}) is an integral of a cohomology element
of $CP^{N-1}\times CP^{N-1}$.

This paper is organized as follows. In Section 2, we introduce the moduli space of lines in $CP^{N-1}$ with one marked point 
and show that it is identified with ${\bf P}({\bf C}^{N}/S_1)$. In Section 3, we introduce the equality (\ref{ingw}) and rewrite 
the r.h.s. of it as an intersection number of ${\bf P}({\bf C}^{N}/S_1)$. In Section 4, we represent the intersection number of   
${\bf P}({\bf C}^{N}/S_1)$ as a residue integral and prove Theorem 1.
\\

\noindent
{\bf Acknowledgment} The author would like to thank Miruko Jinzenji for kind encouragement. 
Research of the author is partially supported by JSPS grant No. 22540061.
   
\section{Moduli Space of Lines in $CP^{N-1}$ with One Marked Point}
First, we introduce a line $l$ in $CP^{N-1}$.
\begin{eqnarray}
&&l:CP^{1}\rightarrow CP^{N-1}\no\\
&&l(s:t)=[{\bf a}_1s+{\bf a}_2t],\;\; 
\label{line} 
\end{eqnarray}
where ${\bf a}_1$ and ${\bf a}_2$ are two linearly independent vectors in ${\bf C}^{N}$.
In (\ref{line}), $[*]$ is used to represent equivalence class under projective equivalence. In this section, 
we consider a line with one marked point $(1:0)\in CP^{1}$, and we introduce $PSL(2,C)$ transformation $\varphi$ that fixes the 
marked point: 
\begin{eqnarray}
&&\varphi: CP^1\rightarrow CP^1\no\\
&&\varphi(s:t)=(\alpha s+\beta t: \delta s),\;\; (\alpha\delta\neq 0) .
\label{tr}
\end{eqnarray}
The moduli space of lines in $CP^{N-1}$ with one marked point is given by the set of equivalence classes of $l$ 
under the equivalence relation:
\begin{equation}
l\circ\varphi\sim l \Longleftrightarrow ({\bf a}_1,{\bf a}_2)\sim (\alpha{\bf a}_1,\beta{\bf a}_1+\delta{\bf a}_2).
\label{mod}
\end{equation} 
From (\ref{mod}), we can see that the moduli space is identified with the projective bundle 
${\bf P}({\bf C}^N/S_1)\rightarrow CP^{N-1}$
where $S_1$ is the tautological line bundle of $CP^{N-1}=\{[{\bf a}_1] \}$. 
From now on, we simply denote by ${\bf P}({\bf C}^N/S_1)$ this projective bundle.
Let $S_2$ be the tautological line bundle of ${\bf P}({\bf C}^N/S_1)$. It is well-known that $H^{*}({\bf P}({\bf C}^N/S_1))$ 
is generated by,
\begin{equation}
h_1:=-c_1(S_1)=c_1(S_1^{*}), \;\;h_2:=-c_1(S_2)=c_1(S_2^{*}).
\end{equation} 
Since the total Chern class of ${\bf C}^{N}/S_1$ is given by,
\begin{equation}
c({\bf C}^{N}/S_1)=\frac{1}{1-th_1}=1+th_1+t^{2}(h_1)^2+\cdots+t^{N-1}(h_1)^{N-1},
\end{equation}
$h_1$ and $h_2$ satisfy the following relation:
\begin{equation}
(h_1)^{N}=0,\;\;\sum_{j=0}^{N-1}(h_1)^{j}\cdot(h_2)^{N-1-j}=0.
\end{equation}
The above relations completely determine the ring structure of $H^{*}({\bf P}({\bf C}^N/S_1))$.
\section{Computation of Number of Lines}
In this section, we denote by $G(2,N)$ the Grassmann variety that parametrizes 2-dimensional linear subspaces in ${\bf C}^{N}$. 
By taking the projective equivalence, it is identified with the moduli space of lines in $CP^{N-1}$ with no marked points. 
$G(2,N)$ has the Scubert cycle $\sigma_a$ (a: non-negative integer) that is defined by 
\begin{eqnarray}
\sigma_{a}:=\{\Lambda\in G(2,N)\;\;|\;\; \dim(\Lambda\cap V_{N-1-a})\geq 1\},
\end{eqnarray} 
where $V_{N-1-a}$ is a fixed $(N-1-a)$-dimensional linear subspace in ${\bf C}^{N}$. 
By taking the projective equivalence, it consists of lines that intersect with the fixed codimension $(a+1)$ linear subspace in $CP^{N-1}$. 
$\sigma_a$ is the homology cycle of 
$G(2,N)$ of complex codimension $a$, and we denote by $\sigma_a^{*}\in H^{*}(G(2,N))$ the Poincare dual of $\sigma_a$. Let $S_G$ be
the tautological rank $2$ vector bundle on $G(2,N)$. We have the following natural projection: 
\begin{eqnarray}
&&\eta:{\bf P}({\bf C}^{N}/S_1)\rightarrow G(2,N),\no\\
&&\eta([({\bf a}_1,{\bf a}_2)])=<{\bf a}_1,{\bf a}_2>_{\bf C}.
\label{projection}
\end{eqnarray}
Here $[({\bf a}_1,{\bf a}_2)]$ is the equivalence class of $({\bf a}_1,{\bf a}_2)$ under (\ref{mod}).
This projection map is equivalent to changing the equivalence relation of pairs of linear independent vectors $({\bf a}_1,{\bf a}_2)$ from 
\begin{equation}
({\bf a}_1,{\bf a}_2)\sim (\alpha{\bf a}_1,\beta{\bf a}_1+\delta{\bf a}_2),\;\;(\alpha\delta\neq0),
\label{sub}
\end{equation} 
into,
\begin{equation}
({\bf a}_1,{\bf a}_2)\sim (\alpha{\bf a}_1+\gamma{\bf a}_2,\beta{\bf a}_1+\delta{\bf a}_2),\;\;(\alpha\delta-\beta\gamma\neq 0).
\label{psl2}
\end{equation}
The equivalence relation (\ref{sub}) is nothing but taking quotient by the subgroup $B_2$ of $PSL(2,{\bf C})(=\mbox{Aut}(CP^1))$ that fixes the marked point $(1:0)$. 
Therefore, the fiber of the map $\eta$ is given by $PSL(2,C)/B_2=CP^1$ and it corresponds to the position of $l(1:0)$ in $l(CP^1)$.

Since $\langle\prod_{j=1}^{n}{\cal O}_{h^{a_j}}\rangle_{0,1}$ is the number of lines in $M_{N}^{k}$ that intersect with codimension 
$a_{j}$ linear subspace in $CP^{N-1}\;\; (j=1,\cdots,n)$, it can be represented as the following intersection number of $G(2,N)$: 
\begin{prop}
\begin{equation}
\langle\prod_{j=1}^{n}{\cal O}_{h^{a_j}}\rangle_{0,1}=
\int_{G(2,N)}c_{top}(S^k(S_G^{*}))\cdot(\prod_{j=1}^{n}\sigma^{*}_{a_j-1}\bigr),
\label{propgw}
\end{equation}
where $S^k(S_G^{*})$ is the $k$-th symmetric product of the dual vector bundle of $S_G$.
\end{prop}
On the other hand, $(h_1)^a \in H^{*}({\bf P}({\bf C}^N/S_1))$ imposes the condition that $l(1:0)\in CP^{N-1}$ lies inside the codimension $a$ linear subspace
of $CP^{N-1}$ and the fiber of $\eta$ is the position of $l(1:0)$ in $l(CP^{1})$. Therefore, we have 
\begin{equation}
\eta_{*}((h_1)^{a})=\sigma_{a-1}^{*}
\label{fi}
\end{equation}
where $\eta_{*}:H^{*}({\bf P}({\bf C}^N/S_1))\rightarrow H^{*}(G(2,N))$ is the fiber integration. Therefore, we obtain from 
the projection formula,
\begin{prop}
\begin{equation}
\langle\prod_{j=1}^{n}{\cal O}_{h^{a_j}}\rangle_{0,1}=
\int_{{\bf P}({\bf C}^N/S_1) }\eta^{*}(c_{top}(S^k(S_G^{*})))\cdot(h_1)^{a_1}\cdot\bigl(\prod_{j=2}^{n}\eta^{*}(\sigma^{*}_{a_j-1})\bigr).
\label{propgw2}
\end{equation}
\end{prop}
\begin{prop}
\begin{eqnarray}
\eta^{*}(c_{top}(S^k(S_G^{*})))&=&\prod_{j=0}^{k}(jh_1+(k-j)h_2)=e^{k}(h_1,h_2),\no\\
\eta^{*}(\sigma_a^{*})&=&\sum_{j=0}^{a}(h_1)^j\cdot (h_2)^{a-j}=w_{a+1}(h_1,h_2).
\label{prope}
\end{eqnarray}
\end{prop}
{\it proof)}
From the definition of $\eta$, we can easily see that, 
\begin{equation}
\eta^{-1}(S_{G})=S_1\oplus S_2.
\label{etai}
\end{equation}
Therefore, the first equality in (\ref{prope}) follows from $\eta^{*}(c(S_G^{*}))=c(S_1^{*})c(S_2^{*})=(1+th_1)(1+th_2)$.
On the other hand, we have the following identity that can be found in the page 411 of \cite{gh}:
\begin{equation}
c({\bf C}^{N}/S_G)=\sum_{a=0}^{N-2}t^{a}\sigma_{a}^{*}.
\label{sigma}
\end{equation}
From these two equalities, we have,
\begin{equation}
\eta^{*}(c({\bf C}^{N}/S_G))=\sum_{a=0}^{N-2}t^{a}\eta^{*}(\sigma_{a}^{*})=\frac{1}{c(S_1)c(S_2)}=\frac{1}{(1-th_1)(1-th_2)}.
\end{equation}
(\ref{prope}) directly follows from this equality.$\Box$
\section{Proof of Theorem 1}
\begin{prop}
Let $f(h_1,h_2)\in H^{*}({\bf P}({\bf C}^N/S_1))$ given as a polynomial in $h_1,h_2$. Then we have,
\begin{equation}
\int_{{\bf P}({\bf C}^N/S_1)}f(h_1,h_2)=\frac{1}{(2\pi\sqrt{-1})^2}\oint_{C_0}dz_2
\oint_{C_0}dz_1\frac{z_2-z_1}{(z_{1})^N\cdot((z_{2})^N-(z_1)^N)}\cdot f(z_1,z_2),
\label{res}
\end{equation}
where $\displaystyle{\frac{1}{2\pi\sqrt{-1}}\oint_{C_0}dz}$ represents the operation of taking a residue at $z=0$.
\end{prop}
{\it proof)}
We take integration paths of $z_1$ and $z_2$ as follows:
\begin{equation}
\{z_1\in{\bf C}\;|\;|z_1|=\epsilon_1\;\},\;\;\{z_2\in{\bf C}\;|\;|z_2|=\epsilon_2\;\}\;\;(0<\epsilon_1<\epsilon_2).
\label{path}
\end{equation}
Under the above condition, $z_1\neq z_2$ and we have,
\begin{equation}
\sum_{j=0}^{N-1}(z_1)^{j}(z_2)^{N-1-j}=\frac{(z_2)^{N}-(z_1)^{N}}{z_2-z_1}.
\label{ratio}
\end{equation}
From degree counting, the r.h.s. of (\ref{res}) is non-zero only if $f(z_1,z_2)$ is a homogeneous polynomial 
of degree $2N-3$ in $z_1$ and $z_2$. Even when it is a homogeneous polynomial of degree $2N-3$, the r.h.s. of 
(\ref{res}) vanishes if it is divided by $(z_1)^N$ or $\displaystyle{\frac{(z_2)^{N}-(z_1)^{N}}{z_2-z_1}}$.
Therefore, it suffices for us to check,
\begin{equation}
\frac{1}{(2\pi\sqrt{-1})^2}\oint_{C_0}dz_2
\oint_{C_0}dz_1\frac{z_2-z_1}{(z_{1})^N\cdot((z_{2})^N-(z_1)^N)}\cdot(z_1)^{N-1}(z_2)^{N-2}=1.
\label{res2}
\end{equation}
But it is obvious because we first integrate $z_1$-variable. $\Box$\\
With this set-up, Theorem 1 follows from Proposition 2, 3, 4 and the following equality:
\begin{prop}
Let $f(z_1, z_2)$ be a polynomial in $z_1$ and $z_2$. Then we have,
\begin{eqnarray}
&&\frac{1}{(2\pi\sqrt{-1})^2}\oint_{C_0}dz_2
\oint_{C_0}dz_1\frac{z_2-z_1}{(z_{1})^N\cdot((z_{2})^N-(z_1)^N)}\cdot f(z_1,z_2)\no\\
&&=\frac{1}{(2\pi\sqrt{-1})^2}\oint_{C_0}dz_2
\oint_{C_0}dz_1\frac{z_2-z_1}{(z_{1})^N\cdot(z_{2})^N}\cdot f(z_1,z_2),
\label{fin}
\end{eqnarray}
where integration paths are taken as given in (\ref{path}).
\end{prop}
{\it proof)}
Let us consider the following integral:
\begin{equation}
\frac{1}{(2\pi\sqrt{-1})^2}\oint_{C_0}dz_2
\oint_{C_0}dz_1\frac{z_2-z_1}{(z_{1})^N\cdot((z_{2})^N-u(z_1)^N)}\cdot f(z_1,z_2),\;\;(0\leq u\leq 1).\no\\
\end{equation}
It suffices for us to show that the above integral is independent of $u$. Under the condition (\ref{path}), 
$\displaystyle{|\frac{z_{1}}{z_{2}}|<1}$ and we can expand the integrand around $z_1=0$:
\begin{equation}
\frac{z_2-z_1}{(z_{1})^N\cdot((z_{2})^N-u(z_1)^N)}\cdot f(z_1,z_2)
=\frac{z_2-z_1}{(z_{1})^N\cdot(z_{2})^N}\cdot f(z_1,z_2)\cdot\biggl(\sum_{j=0}^{\infty}u^j\bigl(\frac{z_1}{z_2}\bigr)^{Nj}\biggr).
\label{texp}
\end{equation}
But the terms with positive power of $u$ have zero contributions to the residue at $z_1=0$ because it contains no negative power of $z_1$. $\Box$

\end{document}